\documentclass[12pt]{amsart}
\usepackage{amscd,amssymb}

\usepackage[graph,frame,poly,arc]{xy}  
\usepackage[plainpages,backref,urlcolor=blue]{hyperref}

\theoremstyle{plain}
\newtheorem{thm}[subsection]{Theorem}

\newtheorem{prop}[subsection]{Proposition}
\newtheorem{cor}[subsection]{Corollary}

\theoremstyle{definition}
\newtheorem{rk}[subsection]{Remark}
\newtheorem{definition}[subsection]{Definition}
\newtheorem{ex}[subsection]{Example}

\numberwithin{equation}{section}
\newcommand{\OO}{{\mathcal O}}

\newcommand{\A}{{\mathcal A}}
\newcommand{\B}{{\mathcal B}}

\newcommand{\al}{{\alpha}}

\newcommand{\Z}{\mathbb{Z}}

\newcommand{\C}{\mathbb{C}}
\newcommand{\PP}{\mathbb{P}}

\DeclareMathOperator{\reg}{reg}
\DeclareMathOperator{\indeg}{indeg}

\begin{document}

\title [On type three complex plane curves]
{On type three complex plane curves}

\author[Alexandru Dimca]{Alexandru Dimca$^1$}
\address{Universit\'e C\^ ote d'Azur, CNRS, LJAD, France and Simion Stoilow Institute of Mathematics,
P.O. Box 1-764, RO-014700 Bucharest, Romania}
\email{Alexandru.Dimca@univ-cotedazur.fr}

\author[Gabriel Sticlaru]{Gabriel Sticlaru}
\address{Faculty of Mathematics and Informatics,
Ovidius University
Bd. Mamaia 124, 900527 Constanta,
Romania}
\email{gabriel.sticlaru@gmail.com }

\thanks{$^1$ partial support from the project ``Singularities and Applications'' - CF 132/31.07.2023 funded by the European Union - NextGenerationEU - through Romania's National Recovery and Resilience Plan.}

\subjclass[2010]{Primary 14H50; Secondary  14B05, 13D02, 32S22}

\keywords{Jacobian ideal, Jacobian module, exponents,  line arrangements, low degree curves}

\begin{abstract} 
The type of a complex projective plane curve has been recently introduced by T. Abe, P. Pokora and the first author. In the same paper they have studied the type two curves. In this paper we study plane curves of type three, with special attention to low degree curves and line arrangements.
\end{abstract}
 
\maketitle

\section{Introduction} 

Let $S=\C[x,y,z]$ be the polynomial ring in three variables $x,y,z$ with complex coefficients, and let $C:f=0$ be a reduced curve of degree $d\geq 3$ in the complex projective plane $\PP^2$. 
We denote by $J_f$ the Jacobian ideal of $f$, i.e. the homogeneous ideal in $S$ spanned by the partial derivatives $f_x,f_y,f_z$ of $f$, and  by $M(f)=S/J_f$ the corresponding graded quotient ring, called the Jacobian (or Milnor) algebra of $f$.
Consider the graded $S$-module of Jacobian syzygies of $f$, namely
$$D_0(f)=\{(a,b,c) \in S^3 \ : \ af_x+bf_y+cf_z=0\}.$$

We say that $C:f=0$ is an {\it $m$-syzygy curve} if  the module $D_0(f)$ is minimally generated by $m$ homogeneous syzygies, say $r_1,r_2,...,r_m$, of degrees $d_j=\deg r_j$ ordered such that $$1 \leq d_1\leq d_2 \leq ...\leq d_m.$$ 
We call these degrees as the {\it exponents} of the curve $C$, and $r_1,...,r_m$ as a {\it minimal set of generators } for the module  $D_0(f)$.

Let $I_f$ denote the saturation of the Jacobian ideal $J_f$ with respect to the maximal ideal ${\bf m}=(x,y,z)$ in $S$. Consider the following local cohomology group, usually called the Jacobian module of $f$, 
 $$N(f)=I_f/J_f=H^0_{\bf m}(M(f)).$$
We set $n(f)_k=\dim N(f)_k$ for any integer $k$ and we introduce the {\it freeness defect of the curve} $C$ by the formula
$$\nu(C)=\max _j \{n(f)_j\}$$ as in \cite{3syz}.
Note that $C$ is free, that is $D_0(f)$ is a free $S$-module of rank 2, if and only if $N(f)=0$ and hence $\nu(C)=0$.
Recall that the total Tjurina number of a given curve $C \subset \mathbb{P}^{2}$  is defined to be
$$\tau(C) = \sum_{p \in {\rm Sing}(C)} \tau_{p},$$ 
where ${\rm Sing}(C)$ denotes the set of all singular points of $C$ and
$\tau_p$  is the Tjurina number of the singularity $(C,p)$.

Consider the general form of the minimal resolution for the Milnor algebra $M(f)$ of a curve $C:f=0$ that is assumed to be not free, namely
\begin{equation}
\label{res2A}
0 \to \oplus_{i=1} ^{m-2}S(-c_i) \to \oplus_{j=1} ^mS(1-d-d_j)\to S^3(1-d)  \to S,
\end{equation}
with $c_1\leq ...\leq c_{m-2}$ and $d_1\leq ...\leq d_m$.
It follows from \cite[Lemma 1.1]{HS} that one has
\begin{equation}
\label{res2B}
c_j=d+d_{j+2}-1+\epsilon_j,
\end{equation}
for $j=1,...,m-2$ and some integers $\epsilon_j\geq 1$. Using \cite[Formula (13)]{HS}, it follows that one has
\begin{equation}
\label{res2C}
d_1+d_2=d-1+\sum_{i=1} ^{m-2}\epsilon_j.
\end{equation}

The following invariant, which leads to an interesting hierarchy among all plane curves, was introduced in \cite{ADP}.

\begin{definition}
\label{deftypek}
Let $C:f=0$ be a reduced plane curve of degree $d$ and let $d_1$ and $d_2$ be the two smallest degrees of a minimal system of generators for the module of Jacobian syzygies $D_0(f)$ as above.
Then we say that $C$ has type $t(C)$, where the positive integer $t(C)$ is defined by the formula
$$t(C)=d_1+d_2+1-d.$$
\end{definition}

The fact that this invariant plays a key role in the study of plane curves
follows from the following alternative definition of it. Let $C:f=0$ be a reduced curve in $\PP^2$ and $r_1$ be a minimal degree generator of the graded $S$-module $D_0(f)$. Then the corresponding Bourbaki ideal $B(C,r_1) \subset S$ 
of $D_0(f)$ was introduced in 
\cite[Section 5]{Split}, and \cite [Section 3]{3syz}, after being considered
 in a non-explicit way by du Plessis–Wall in \cite{duPCTC} and the first author in \cite{Dmax}. Its properties were further explored in the recent paper \cite{JNS}. It follows from the definition of the Bourbaki ideal $B(C,r_1)$ that the type of the curve $C$ is just the initial degree of the
homogeneous ideal $B(C,r_1)$, namely
\begin{equation}
\label{eqT}
t(C)=\indeg B(C,r_1)= \min \{s \in \Z \ : \ B(C,r_1)_s \ne 0 \}.
\end{equation}

Note that a curve $C$ is free exactly when $t(C)=0$, and a curve $C$ is plus-one generated  as defined in \cite{Abe, 3syz} exactly when $t(C)=1$. 
The curves $C$ such that $t(C)=2$ have been studied in \cite{ADP}.

In this note we consider the curves satisfying $t(C)=3$. In Section 2 we recall some basic facts  related to Jacobian syzygies of curves and in Section 3 we use these facts to derive formulas for the Tjurina numbers, see Theorem \ref{thm1}, and freeness defects, see Propositions \ref{thm3A}, \ref{thm3B}, \ref{thm3B'} and \ref{thm3C}, for curves of type 3. In particular we show that these curves fall naturally into 4 subclasses, denoted by $3A$, $3B$, $3B'$ and $3C$. The fact that all these 4 subclasses do exist is shown by Example \ref{exABC}.

In Section 4 we relate the curves of type 3 to other classes of curves, namely the maximal Tjurina curves, see Proposition \ref{propMAX} and Example \ref{exC3}, and to the nodal curves, to which the longer subsection (4.4) is devoted. In particular, many nodal curves of degree 6 and type 3 are constructed in Example \ref{exNodal}.
Then we show that the Zariski sextic with 6 cusps on a conic and
Ziegler-Yuzvinski arrangement of 9 lines having 6 triple points on a conic are both curves of type 3, see Examples \ref{exOC1} and \ref{exZ}. At the end of this section we construct a countable family of curves of type $3C$ enjoying nice symmetry properties, see Theorem \ref{thmSym}.

In Section 5 we classify up to a certain point the curves of degree 4 (resp. degree 5) having type 3, see Proposition \ref{propd4} and Corollary \ref{cord4} (resp. Proposition \ref{propd5} and Corollary \ref{cord5}).

In Section 6 we consider the line arrangements of type 3. First, in Propositions \ref{propA6}, \ref{propA7} and resp. \ref{propA8} we classify up to a certain point the  arrangements of 6, 7 and resp. 8 lines that have type 3. Finally, in Theorem \ref{thmAd} we construct a arrangement of $d=n_1+n_2+2$ lines and having type $3C$, for any integers $2 \leq n_1 \leq n_2$. 

\section{Basic facts  related to Jacobian syzygies of curves }

We recall first the inequalities
\begin{equation} 
\label{eqd-1}
d_1 \leq d_2 \leq d_3 \leq d-1
\end{equation} 
which hold for any reduced plane curve, see \cite[Theorem 2.4]{3syz}.
If the generator $r_j$ of $D_0(f)$ has components $(a_j,b_j,c_j)$ then one has
$$a_j,b_j,c_j \in S_{d_j} \text{ and } a_jf_x+b_jf_y+c_jf_z=0,$$
for any $j=1, \ldots, m$, where $(d_1,d_2, \ldots,d_m)$ are the exponents of $C:f=0$. Moreover, the second order syzygies are given by equalities in $S^3$ of the form
\begin{equation} 
\label{syz2}
\al_{i,1}r_1+ \al_{i,2}r_2+ \ldots + \al_{i,m}r_m=0
\end{equation} 
with $i=1, \ldots ,m-2$,  and $\al_{i,j}$ homogeneous polynomials in $S$ such that if $\al_{i,j}\ne 0$ then
\begin{equation} 
\label{syz21}
\deg(\al_{i,j})=d_m-d_j+\epsilon_i,
\end{equation} 
for any $i=1, \ldots ,m-2$ and any $j=1, \ldots, m$.

We recall next the following result due to du Plessis and Wall, see \cite[Theorem 3.2]{duPCTC} as well as \cite{E} for an alternative approach.
\begin{thm}
\label{thmCTC}
For positive integers $d$ and $d_1$, define two new integers by 
$$\tau(d,d_1)_{min}=(d-1)(d-d_1-1)  \text{ and } 
\tau(d,d_1)_{max}= (d-1)^2-d_1(d-d_1-1).$$ 
If $C:f=0$ is a reduced curve of degree $d$ in $\PP^2$ and  $d_1$ is the first exponent of $f$,  then
$$\tau(d,d_1)_{min} \leq \tau(C) \leq \tau(d,d_1)_{max}.$$
Moreover, for $d_1 \geq d/2$, the stronger inequality
$\tau(C) \leq \tau'(d,d_1)_{max}$ holds, where
$$\tau(d,d_1)'_{max}=\tau(d,d_1)_{max} - \binom{2d_1+2-d}{2}.$$
\end{thm}
The following result describes the curves with a minimal Tjurina number, see
\cite[Theorem 3.5 (1)]{3syz}.
\begin{thm}
\label{thmMT}
If $C:f=0$ is a reduced curve of degree $d$ in $\PP^2$ and  $d_1$ is the first exponent of $f$,  then $\tau(C)=\tau(d,d_1)_{min}$ if and only if
$C$ is a 3-syzygy curve with $d_2=d_3=d-1$.
\end{thm}

With the notation and assumptions from \eqref{res2A}, the Castelnuovo-Mumford regularity $reg(D_0(f))$ of the $S$-module $D_0(f)$ is given by 
\begin{equation}
\label{eqCM}
\reg D_0(f)=d_m+\epsilon_{m-2}-1.
\end{equation}
This invariant satisfies
\begin{equation}
\label{eqCM2}
\reg D_0(f)\leq 2d-4,
\end{equation}
for $d \geq 3$, see \cite[Remark 5.1]{CMreg}. 
In view of \eqref{eqCM}, this implies in particular
\begin{equation}
\label{eqCM3}
d_m \leq 2d-4,
\end{equation}
In addition, equality holds in \eqref{eqCM2} and \eqref{eqCM3} when $C$ is a curve having only one node as singularities, that is when $\tau(C)=1$.
Moreover, we have the following result, \cite[Theorem 2.1]{CMreg}.
\begin{thm}
\label{thmCM}
Assume that $C$ has $s$ irreducible components $C_i$, of degree $e_i$, which are all smooth, for $i=1,\ldots,s$. Let $e$ be the maximum of the degrees $e_i$. Then
$$\reg D_0(f) \leq d+e-3$$
and equality holds when $C$ is a nodal curve.

\end{thm}

Finally, we recall the following result, see \cite[Theorem 1.2]{Dimca1}.

\begin{thm}
\label{deff}
Let $C \, : f=0$ be a reduced plane curve of degree $d$ and $d_{1}= {\rm mdr}(C)$. Then the following hold.
\begin{enumerate}
    \item If $d_{1} \leq (d-1)/2$, then $\nu(C) = (d-1)^{2}-d_{1}(d-1-d_{1})-\tau(C)$.
    \item If $d_{1} \geq (d-1)/2$, then
    $$\nu(C) = \bigg\lceil\frac{3}{4}(d-1)^{2}\bigg\rceil - \tau(C).$$
\end{enumerate}
\end{thm}

\section{Tjurina numbers and freeness defects} 

The exponents and the corresponding Tjurina numbers of the curves of type $3$ are described in the following result, where the notation from \eqref{res2C} is used.

\begin{thm}
\label{thm1}
Let $C$ be a curve of type $3$.
Then one of the following four  cases occurs.

\begin{enumerate}

\item $C$ is a 3-syzygy curve and $\epsilon_1=3$. If $(d_1,d_2,d_3)$ are the exponents of $C$, then
$$\tau(C)=d_1^2+d_1d_2+d_2^2-3(d_1+d_2+d_3).$$

\item $C$ is a 4-syzygy curve and  $\epsilon_1= 2, \epsilon_2= 1$.
If $(d_1,d_2,d_3,d_4)$ are the exponents of $C$, then $d_3<d_4$ and
$$\tau(C)=d_1^2+d_1d_2+d_2^2-3d_1-3d_2-2d_3-d_4+2.$$

\item $C$ is a 4-syzygy curve and  $\epsilon_1= 1, \epsilon_2= 2$.
If $(d_1,d_2,d_3,d_4)$ are the exponents of $C$, then
$$\tau(C)=d_1^2+d_1d_2+d_2^2-3d_1-3d_2-d_3-2d_4+2.$$

\item $C$ is a 5-syzygy curve and  $\epsilon_1=  \epsilon_2=  \epsilon_3=1$. 
If $(d_1,d_2,d_3,d_4,d_5)$ are the exponents of $C$, then
$$\tau(C)=d_1^2+d_1d_2+d_2^2-3d_1-3d_2-d_3-d_4-d_5+3.$$

\end{enumerate}
Conversely, if  one of the following occurs

\begin{enumerate}

\item[(a)] $C$ is a 3-syzygy curve and $\epsilon_1=3$,

\item[(b)]  $C$ is a 4-syzygy curve and  $\epsilon_1= 2, \epsilon_2= 1$,

\item[(c)]  $C$ is a 4-syzygy curve and  $\epsilon_1= 1, \epsilon_2= 2$,

\item[(d)]  $C$ is a 5-syzygy curve and  $\epsilon_1=  \epsilon_2=  \epsilon_3=1$,

\end{enumerate}
 then $t(C)=3$.

\end{thm}
In case $(1)$ (resp. $(2)$, $(3)$, $(4)$), {\it we say that $C$ is a curve of type} $3A$ (resp. $3B$, $3B'$, $3C$).
The additional condition $d_3<d_4$ in case (2) comes from the inequality 
$c_1 \leq c_2$ in \eqref{res2A}, which is equivalent to
$$d_3+2=d_3+\epsilon_1 \leq d_4 + \epsilon_2=d_4+1.$$

\proof

To see that the number $m$ of syzygies satisfies $3 \leq m \leq 5$, and to get the corresponding values for the $\epsilon_j$, use the  equality \eqref{res2C}. In case (1), the value for $\tau(C)$ follows from
\cite[Proposition 2.1]{3syz} using that $d_1+d_2=d+2$, where $d= \deg C$. 

In the  case  (2), we use \eqref{res2B} to get the following resolution of the Milnor algebra $M(f)$:
$$0 \to S(-d-d_3-1)\oplus S(-d-d_4)  \to \oplus_{j=1}^{5}S(1-d-d_j) \to S^3(1-d) \to S.$$
This resolution implies that in this case one has
$$\tau(C)=\binom{k+2}{2}-3 \binom{k+3-d}{2}+\sum_{j=1}^{4}\binom{k+3-d-d_j}{2}-$$
$$-  \binom{k+1-d-d_3}{2}-  \binom{k+2-d-d_4}{2}$$
for any large $k$.
In the  case  (3), we use \eqref{res2B} to get the following resolution of the Milnor algebra $M(f)$:
$$0 \to S(-d-d_3)\oplus S(-d-d_4-1)  \to \oplus_{j=1}^{5}S(1-d-d_j) \to S^3(1-d) \to S.$$
This resolution implies that in this case one has
$$\tau(C)=\binom{k+2}{2}-3 \binom{k+3-d}{2}+\sum_{j=1}^{4}\binom{k+3-d-d_j}{2}-$$
$$-  \binom{k+2-d-d_3}{2}-  \binom{k+1-d-d_4}{2}$$
for any large $k$.

In case (4), we use \eqref{res2B} to get the following resolution of the Milnor algebra $M(f)$:
$$0 \to S(-d-d_3)\oplus S(-d-d_4) \oplus S(-d-d_5) \to \oplus_{j=1}^{5}S(1-d-d_j) \to S^3(1-d) \to S.$$
This resolution implies that in this case one has
$$\tau(C)=\binom{k+2}{2}-3 \binom{k+3-d}{2}+\sum_{j=1}^{5}\binom{k+3-d-d_j}{2}-$$
$$-  \binom{k+2-d-d_3}{2}-  \binom{k+2-d-d_4}{2}-  \binom{k+2-d-d_5}{2}$$
for any large $k$.
The last claim follows from the equality \eqref{res2C}.
\endproof

For the curves $C$ of type 3, we also have simple formulas for the freeness defect $\nu(C)$, as shown below.

\begin{prop}
\label{thm3A}
Suppose that $C$ has type $3A$ and  exponents $(d_1,d_2,d_3)$.
Then the following hold.
\begin{enumerate}

\item If $d_1 \leq d_2-2$, then
$\nu(C)=3(d_3-d_2)+9 \geq 9.$

\item If $d_1=d_2-1$, then
$\nu(C)= 3(d_3-d_2)+8 \geq 8.$
\item If $d_1=d_2$, then
$\nu(C)= 3(d_3-d_2)+7\geq 7.$

\end{enumerate}

\end{prop}

\proof

We  use Theorem \ref{deff} and Theorem \ref{thm1}. Note that the condition $d_1 \leq d_2-3$ is exactly the condition $d_1\leq (d-1)/2$ which occurs in Theorem \ref{deff} (1). Hence, to prove the claim $(1)$ in Proposition \ref{thm3A} under this stronger condition, we use the formula for $\nu(C)$ given in Theorem \ref{deff} (1) and the formula for $\tau(C)$ given in  Theorem \ref{thm1} (1). To prove the remaining part of claim (1) as well as the claims (2), (3) and (4), we use 
the formula for $\nu(C)$ given in
Theorem \ref{deff} (2) and the formula for $\tau(C)$ given in Theorem \ref{thm1} (1).  This completes the proof of Proposition \ref{thm3A}.

\endproof

The proofs of  Propositions \ref{thm3B},  \ref{thm3B'} and  \ref{thm3C}  below are analogous to the above proof and we leave them to the reader.

\begin{prop}
\label{thm3B}
Suppose that $C$ has type $3B$ and  exponents $(d_1,d_2,d_3,d_4)$.
Then the following hold.
\begin{enumerate}

\item If $d_1 \leq d_2-2$, then
$\nu(C)=2d_3+d_4-3d_2+7 \geq 7.$

\item If $d_1=d_2-1$, then
$\nu(C)= 2d_3+d_4-3d_2+6 \geq 6.$
\item If $d_1=d_2$, then
$\nu(C)= 2d_3+d_4-3d_2+ 5 \geq 5.$

\end{enumerate}

\end{prop}
\begin{prop}
\label{thm3B'}
Suppose that $C$ has type $3B'$ and  exponents $(d_1,d_2,d_3,d_4)$.
Then the following hold.
\begin{enumerate}

\item If $d_1 \leq d_2-2$, then
$\nu(C)=d_3+2d_4-3d_2+7 \geq 7.$

\item If $d_1=d_2-1$, then
$\nu(C)= d_3+2d_4-3d_2+6 \geq 6.$
\item If $d_1=d_2$, then
$\nu(C)= d_3+2d_4-3d_2+5\geq 5.$

\end{enumerate}

\end{prop}
\begin{prop}
\label{thm3C}
Suppose that $C$ has type $3C$ and  exponents $(d_1,d_2,d_3,d_4,d_5)$.
Then the following hold.
\begin{enumerate}

\item If $d_1 \leq d_2-2$, then
$\nu(C)=d_3+d_4+d_5-3d_2+6 \geq 6.$
\item If $d_1=d_2-1$, then
$\nu(C)= d_3+d_4+d_5-3d_2+5 \geq 5.$
\item If $d_1=d_2$, then
$\nu(C)= d_3+d_4+d_5-3d_2+4\geq 4.$

\end{enumerate}

\end{prop}

Let us look at the first examples of curves $C$ having type $3$.

\begin{ex}
\label{exABC}
(i) The curve $C: f=xy(x^4 + y^4-z^4)=0$ has type $3A$ and exponents $(4,4,5)$, as a direct computation using for instance SINGULAR \cite{Singular} shows. Using Theorem \ref{thm1} (1) we see that $\tau(C)=9$ and Proposition \ref{thm3A} (3) implies that $\nu(C)=10$.

(ii) The curve $C: f=x^4+y^4+xyz(x+z)=0$ has type $3B$ and exponents $(3,3,3,4)$. Using Theorem \ref{thm1} (2) we see that $\tau(C)=1$ and Proposition \ref{thm3B} (3) implies that $\nu(C)=6$.

(iii) The curve $C: f=(x^2+y^2-2z^2)(x^9+y^9-z(x^2-3y^2)^4)=0$ has type $3B'$ and exponents $(6,7,8,9)$. Using Theorem \ref{thm1} (3) we see that $\tau(C)=64$ and Proposition \ref{thm3B'} (2) implies that $\nu(C)=11$.

(iv) The line arrangement
$$C: f= (x-2y)(y-3z)(z-5x)(x-y+z)(x+y-z)(x+2y-3z)$$
$$(2x+3y-5z)(3x+5y-8z)(5x+7y-12z)=0$$
 has type $3C$ and exponents $(5,6,7,7,7)$. Using Theorem \ref{thm1} (4) we see that $\tau(C)=40$ and Proposition \ref{thm3C} (2) implies that $\nu(C)=8$.

\end{ex}

\begin{rk}
\label{rkABC}
In the case of curves of type 3 we can say more about the second order syzygies \eqref{syz2}. For instance, if $C:f=0$ has type $3A$, then there is a unique generating second order syzygy since $m=3$, and the equalities \eqref{syz21} yield the following inequalities
$$\deg \al_{1,1} \leq d-1, \  \deg \al_{1,2} \leq \frac{d}{2}+1 \text{ and } \deg \al_{1,3}=3.$$
For the first inequality we use the fact that $d_1 \geq 3$ which follows from \eqref{eqd-1} and $t(C)=3$. For the second inequality we use
$d_2 \geq (d+2)/2$, and for the last equality we use the fact that
$\epsilon_1=3$ in this case as in Theorem \ref{thm1} (1).
Similar relations may be obtained for the curves of type $3B$, $3B'$ and $3C$.
\end{rk}

\section{Geometric examples of curves of type 3} 

\subsection{Maximal Tjurina curves} 
 A reduced curve $C:f=0$
of degree $d$, such that $d_1 \geq d/2$ and the equality
\begin{equation} 
\label{eqMAX}
\tau(C) =  \tau(d,d_1)'_{max},
\end{equation}
with $\tau(d,d_1)'_{max}$ as in Theorem \ref{thmCTC}, holds for $C$ is called a {\it maximal Tjurina curve of type} $(d,d_1)$,
 see  \cite{maxTjurina}. Such a curve satisfies $m=2d_1-d+3$ and
 $$d_1=d_2= \ldots = d_m.$$
 It follows that 
 $$t(C)=2d_1-d+1=m-2.$$
 Hence to get examples of curves with $t(C)=3$ in this class of curves we need
 $m=5$. In fact we have the following.
 \begin{prop} 
 \label{propMAX}
 For a reduced curve $C:f=0$
of degree $d$, the following are equivalent.

\begin{enumerate}
    \item  The curve $C$ has type $3C$ and its exponents satisfy $d_1=d_2= \ldots =d_5 $.
    \item The curve $C$ is maximal Tjurina of type $(d,d_1)$ with $2d_1=d+2$.
\end{enumerate} 
 \end{prop}
\proof
Assume first that $C$ satisfies (1). It follows that $2d_1=d_1+d_2=d+2$, and hence $d_1 \geq d/2$.
Then, using Theorem \ref{thm1} (4), we get
$$\tau(C)=3d_1^2-9d_1+3,$$
which shows that the equality \eqref{eqMAX} holds, and hence $C$ is 
maximal Tjurina of type $(d,d_1)$. 

Assume now that $C$ satisfies (2). Then the equality $2d_1=d+2$ implies as explained above that $m=5$ and that $t(C)=3$. Hence $C$ has type $3C$
and moreover $$d_1=d_2= \ldots =d_5.$$
\endproof 
Proposition \ref{propMAX} and the lists of maximal Tjurina curves in \cite{maxTjurina} yield the following examples of curves of type $3C$.

\begin{ex}
\label{exC3}
(i) A quartic $C$ with 3 nodes is a maximal Tjurina curve of type $(4,3)$, see \cite[Proposition 5.8]{maxTjurina}, and hence $C$ is a curve of type $3C$.

(ii) An arrangement $C$ of 6 lines in general position, in other words having only double points, is a maximal Tjurina curve of type $(6,4)$, see \cite[Proposition 5.11]{maxTjurina}, and hence $C$ is a curve of type $3C$.

(iii) Line arrangements which are maximal Tjurina curves of type $(8,5)$ and $(10,6)$ are constructed in  \cite[subsections (5.5) and (5.6)]{maxTjurina}, and hence both such arrangements are  curves of type $3C$.

\end{ex}

\subsection{Nodal curves} 

Recall the following result, see \cite[Theorem 4.1]{Edin}.

\begin{thm}
\label{thmNodal}
Let $C$ be a nodal curve of degree $d$ having $r\geq 2$ irreducible components. Then
$$d_1=\ldots = d_{r-1}=d-2.$$
\end{thm}
We deduce the following direct consequence.
\begin{cor}
\label{corNodal}
Let $C$ be a nodal curve of degree $d$ having $r\geq 3$ irreducible components. Then $C$ has type $3$ if and only if $d=6$.
\end{cor}
\begin{ex}
\label{exNodal}
(i) When $C$ is a nodal arrangement of 6 lines, we have seen in Example
\ref{exC3} (ii) that $C$ has type $3C$.

(ii) The nodal curve
$$C_2: f=xyz(x+y+z)(x^2+y^2+z^2)=0$$
consists of a conic and 4 lines, and a direct computation with SINGULAR shows that $C$ has type $3B'$, namely the exponents are $(4,4,4,4)$.
Similarly, the nodal curve
$$C_3: f=xyz(x^3+y^3+z^3)=0$$
has exponents $(4,4,4)$ and hence has type $3A$,
the nodal curve
$$C_4: f=xy(x^4+y^4+z^4)=0$$
has exponents $(4,4,5)$ and type $3A$. 

On the other hand,  the nodal curve
$$C_5: f=x(x^5+y^5+z^5)=0,$$
with only 2 irreducible components, has exponents $(4,5,5)$ and hence has type 4.

(iii) The nodal curve
$$C_{2,2}: f=xy(x^2+y^2+z^2)(x^2+2y^2+3z^2)=0$$
consists of two smooth conics and two lines, has exponents $(4,4,4,5)$ and
$\tau(C)=13$. It follows that $C_{2,2}$ has type $3B$.
The nodal curve
$$C_{2,3}: f=x(x^2+y^2+z^2)(x^3+2y^3+3z^3)=0$$
consists of a smooth conic, a smooth cubic and a line, has exponents $(4,4,5,5)$ and  hence $C_{2,3}$ has type $3B'$.
 The nodal curve
$$C_{2,2,2}: f=(x^2+y^2+z^2)(x^2+2y^2+3z^2)(x^2-y^2+5z^2)=0$$
consists of 3 smooth conics,  has exponents $(4,4,5,5,5)$ and  hence $C_{2,2,2}$ has type $3C$.

\end{ex}
When $C$ is a nodal curve with at most two irreducible components, we have the following consequence of Theorem \ref{thmNodal}.
\begin{cor}
\label{corNodal2}
Let $C$ be a nodal curve of degree $d$ having one (resp. two) irreducible component (resp. components). Then $C$ has type $3$ if and only if $d=4$ (resp.$d=5$).

\end{cor}
\proof
The claims follow from Theorem \ref{thmNodal} and the inequalities
\eqref{eqd-1}.
\endproof

\begin{ex}
\label{exNodal2}
(i) A smooth or an irreducible nodal quartic has type 3. A smooth quartic has type $3A$ since $m=3$, and an irreducible quartic with 3 nodes has type $3C$ as shown in Example \ref{exC3} (i). The quartic
$$C: f=xyz^2+x^4+y^4=0$$
and the quartic $C$ in Example \ref{exABC} (ii) have a unique node, exponents $(3,3,3,4)$ and $\tau(C)=1$. Hence these curves $C$ have type $3B$.

(ii) The nodal curve
$$C'_{2,3}: f=(x^2+y^2+z^2)(x^3+2y^3+3z^3)=0$$
consists of a smooth conic and  a smooth cubic, has exponents $(3,4,4,4)$ and  hence $C'_{2,3}$ has type $3B'$.
\end{ex}

\subsection{Other curves of type 3} 

\begin{ex}
\label{exOC1}
The sextic
$$C: f=(x^3+y^3)^2+(y^2+z^2)^3=0$$
has 6 cusps $A_2$ and goes back to the work of Zariski on the fundamental groups of complements of projective plane curves, see for instance \cite[Proposition 4.4.16 and Theorem 6.4.9]{STH}. A direct computation shows that this curve has exponents $(3,5,5,5)$, and hence has type $3B'$.

\begin{ex}
\label{exZ}
When $C:f=0$ is the Ziegler-Yuzvinsky arrangement of 9 lines obtained starting with a hexagon inscribed in a conic, see \cite{Z,Y}, the
corresponding exponents are $(5,6,6,6)$ and hence $C$ is a curve of type $3B'$.
As explicit equation for $C$ one may take
$$xyz(x+y-z)(x-y+z)(2x-2y+z)(2x-y-2z)(2x+y+z)(2x-y-z)=0, $$  
when the conic $Q$ is a union of two lines,
or
$$xy(x-y-z)(x-y+z)(2x+y-2z)(x+3y-3z)(3x+2y+3z)(x+5y+5z)(7x-4y-z)=0,$$
when the conic $Q$ is smooth, see for details \cite{DSZi}.
\end{ex}

\end{ex}
We end this section by describing a countable family of curves of type $3C$ having a nice symmetry.
\begin{thm}
\label{thmSym}
The curve 
$$C_j:f=x^jy^j+y^jz^j+x^jz^j=0$$
for $j\geq 3$ has degree $d=2j$ and exponents $(j+1,j+1,2j-1,2j-1,2j-1)$. In particular any such curve has type $3C$.
\end{thm}
\proof
It is easy to check the following 3 equalities
$$xf_x-yf_y+zf_z=2jx^jz^j, \ -xf_x+yf_y+zf_z=2jy^jz^j \text{ and } xf_x+yf_y-zf_z=2jx^jy^j.$$
By multiplying the first equality by $y^j$ and subtracting the second equality (resp. the third) multiplied by $x^j$ (resp. by $z^j$) we get two syzygies $r_1$ and $r_2$ for $f$ of degree $j+1$. Looking at the formulas for $f_x,f_y$ and $f_z$ it is clear that syzygies of degree $<j+1$ do not exist, and hence $d_1=d_2=j+1$. In particular
$C_j$ has type 3.

To show that the next syzygy $r_3$ may occur only in degree $d_3 =2j-1$ we use the equality in \cite[Theorem 1]{Bull} for $k=2j-3$. This yields the equality
\begin{equation} 
\label{eqSym}
\dim D_0(f)_{2j-2}=\tau(C_j)-\dim S_{2j-3}/I_{f,2j-3}.
\end{equation} 
The curve $C_j$ has 3 singular points located at $p_1=(0:0:1)$, $p_2=(0:1:0)$ and $p_3=(1:0:0)$. At $p_1$, the corresponding singularity is obtained from $f$ by setting $z=1$, and hence 
it is given by
$$g_1=x^j+y^j+x^jy^j=u^j+v^j=0,$$
where $u=x\rho$ and $v=y$ with $\rho^j=1+y^j$ an invertible element in the local ring at $p_1$. It follows that
the local
Tjurina number at $p_1$ is $\tau(C_j,p_1)=(j-1)^2$. The same argument applies to $p_2$ and $p_3$ and hence
$$\tau(C_j)=3(j-1)^2.$$
Next the saturation $I_f$ of the Jacobian ideal $J_f$ consists of polynomials such that their localizations at all the singular points $p_j$ lie in the corresponding local Jacobian ideal $J_{g_j}$.
Note that, with the above notation, $J_{g_1}$ is generated by $u^{j-1}$ and $v^{j-1}$, and also by $x^{j-1}$ and $y^{j-1}$. It follows that $h \in I_{f,2j-3}$ is a sum of monomials, each one of then containing at least two variables among $x,y,z$ with an exponent at least $j-1$. Indeed, if only one such variables occurs, say $x$, then by localizing at $p_3$ we get a contradiction. It follows that $I_{f,2j-3}=0$ and hence
$$\dim S_{2j-3}/I_{f,2j-3}=\dim S_{2j-3}= \binom{2j-1}{2}=(j-1)(2j-1).$$
Using \eqref{eqSym} it follows that
$$\dim D_0(f)_{2j-2}=3(j-1)^2-(j-1)(2j-1)=j^2-3j+2.$$
Since $C_j$ has type 3, only one of the following cases is possible.

\medskip

\noindent {\bf Case 1: $C_j$ has type $3A$}.
Then using the formula for $\tau(C_j)$ given in Theorem \ref{thm1} we get
$$3(j-1)^2=\tau(C_j)=3(j+1)^2-6(j+1)-3d_3,$$
which implies
$d_3=2j-2$. But the two syzygies $r_1$ and $r_2$ of degree $j+1$ span a subspace in 
$\dim D_0(f)_{2j-2} $ of dimension
$$2 \dim S_{j-3}=j^2-3j+2,$$
and hence they span the whole subspace $\dim D_0(f)_{2j-2} $.
Indeed, one has 
$$S_{j-3}r_1 \cap S_{j-3}r_2=0,$$
since the $S$-module $D_0(f)/Sr_1$ is torsion free, see the end of the proof of \cite[Theorem 1.4]{ADP}.
This contradiction shows that this case cannot occur.

\medskip

\noindent {\bf Case 2: $C_j$ has type $3B$ or $3B'$}.
Then, if the exponents are $(j+1,j+1,d_3,d_4)$, then using the total Tjurina numbers as above we get
$$2d_3+d_4=6j-4 \text{ or } d_3+2d_4=6j-4.$$
As above we see that $d_4 \geq d_3 >2j-2$, which leads again to a contradiction.

\medskip

\noindent {\bf Case 3: $C_j$ has type $3C$ }.

Then, if the exponents are $(j+1,j+1,d_3,d_4,d_5)$, then using the total Tjurina numbers as above we get
$$d_3+d_4+d_5=6j-3.$$
As above we see that $d_5 \geq d_4 \geq d_3 >2j-2$, therefore the only solution is 
$$d_5 =d_4 =d_3 =2j-1.$$
This completes the proof of Theorem \ref{thmSym}.

\endproof

\begin{rk}
\label{rkSym}
It is interesting that the curves obtained from the curve $C_j$ above by adding some of the coordinate lines $x,y$ or $z$ have smaller types.
Indeed, one can check by SINGULAR and for small values of $j\geq 3$ the following facts.

\medskip
\noindent (i) The curve 
$$C_j':f=x(x^jy^j+y^jz^j+x^jz^j)=0$$
 has degree $d=2j+1$ and exponents $(j+1,j+1,2j-1)$, hence it has type $2A$.

\medskip
\noindent
(ii) The curve 
$$C_j'':f=xy(x^jy^j+y^jz^j+x^jz^j)=0$$
has degree $d=2j+2$ and exponents $(j+1,j+1,2j)$, 
hence it has type 1, namely it is a plus-one generated  curve, see \cite{Abe,3syz}.

\medskip
\noindent
(iii) The curve 
$$C_j''':f=xyz(x^jy^j+y^jz^j+x^jz^j)=0$$
has degree $d=2j+3$ and exponents $(j+1,j+1)$, hence
it has type 0, namely it is a free curve.

\end{rk}

\section{On the curves of type 3 and low degree} 

In this section we give examples of curves $C$ of type 3, {\it without using the 
defining equation $f=0$ or the SINGULAR package}.
In view of  \eqref{eqd-1}, for a curve $C$ of type 3 one has
$$d+2=d_1+d_2 \leq 2d-2,$$
which yields $d\geq 4$. Hence we start our analysis with curves of degree 4.
\begin{prop}
\label{propd4}
If $C$ is a curve of degree $d=4$ and of type 3, then $d_1=d_2=d_3=3$ and
$\tau(C) \leq 3$. Moreover, only the following cases may occur.
\begin{enumerate}
    \item  $\tau(C)=0$ and exponents $(3,3,3)$.
    \item $\tau(C)=1$ and exponents $(3,3,3,4)$.
    \item $\tau(C)=2$ and exponents $(3,3,3,3)$ or $(3,3,3,3,4)$.
    \item $\tau(C)=3$ and exponents  $(3,3,3,3,3)$.
\end{enumerate}
Conversely, any curve $C$ of degree $d=4$ such that $\tau(C) \leq 2$ has type 3. 
\end{prop}
\proof
If $d=4$ and $t(C)=3$, then  \eqref{eqd-1} implies that $d_1=d_2=d_3=3$.
Since $\tau(4,3)'_{max}=3$, Theorem \ref{thmCTC} implies that 
$$0 \leq \tau(C)\leq 3.$$
Using \eqref{eqCM3} we get $d_m \leq 4$, and hence only the exponents listed in Proposition \ref{propd4} may occur. For each exponent we compute the Tjurina number $\tau(C)$ using Theorem \ref{thm1}.
When  $\tau(C)=1$, exponents $(3,3,3,4)$  and  $(3,3,3,4,4)$ seem possible.
However, it is known that a curve of degree $d$ with a single node has exponents $(d-1,d-1,d-1,2d-4)$, see \cite[Remark 2.5]{3syz}.
The last claim follows from Theorem \ref{thmCTC}. Indeed,
$$\tau(4,2)_{min}=3$$
and hence $\tau(C) \leq 2$ implies $d_1=3$. Then $d_2=3$ by \eqref{eqd-1},
and hence $t(C)=3$.
\endproof

\begin{cor}
\label{cord4}
A quartic curve $C$ with $\tau(C)=3$ is either a curve of type $2A$ with exponents $(2,3,3)$, or a curve of type $3C$ with exponents
$(3,3,3,3,3)$.
\end{cor}
\proof
If $d_1=3$, then $C$ is of type $3C$ as it follows from
Proposition \ref{propd4}.
If $d_1=2$, one uses the equality
$$\tau(4,2)_{min}=3=\tau(C)$$
and Theorem \ref{thmMT} to get $d_2=d_3=3$.
Hence  the exponents are $(2,3,3)$ and  $C$ is a curve of type  $2A$, see \cite[Proposition 1.11]{ADP}.
The case $d_1=1$ is impossible, since $\tau(4,1)_{min}=6$.
\endproof

\begin{ex}
\label{exd4}
(i) The curve
$$C: f=x^4+x^2z^2+y^3z=0$$
has a cusp $A_2$ at the point $(0:0:1)$, exponents $(3,3,3,3)$ and $\tau(C)=2$. The curve
$$C: f=xyz(x+y+z)+(x+y)^4=0$$
has 2 nodes at the points $(1:0:0)$ and $(0:1:0)$, exponents $(3,3,3,3)$ and $\tau(C)=2$.
{\it We do not know if there is a curve $C$ with $d=4$ and exponents $(3,3,3,3,4)$.}

(ii) The nodal curve
$$C: f=x(x^3+y^3+z^3)=0$$
has exponents $(2,3,3)$, $t(C)=2$ and $\tau(C)=3$. This shows that the last claim in Proposition \ref{propd4} cannot be improved.
\end{ex}
Next we consider the case of quintic curves.
\begin{prop}
\label{propd5}
If $C$ is a curve of degree $d=5$ and of type 3, then $d_1=3$ and $d_2=d_3=4$ and
$4\leq \tau(C) \leq 7$. Moreover, only the following cases may occur.
\begin{enumerate}
    \item  $C$ has type $3A$, $\tau(C)=4$ and exponents $(3,4,4)$.
    \item $C$ has type $3B$, $\tau(C)=5$ and exponents $(3,4,4,5)$.
    \item $C$ has type $3B'$, $\tau(C)=6$ and exponents $(3,4,4,4)$.
    \item $C$ has type $3C$, and one of the following holds.
    
    $\tau(C)=5$ and exponents $(3,4,4,5,5)$ or $(3,4,4,4,6)$;
    
    $\tau(C)=6$ and exponents $(3,4,4,4,5)$;
    
    $\tau(C)=7$ and  exponents $(3,4,4,4,4)$.

\end{enumerate}
Conversely, any nodal curve $C$ of degree $d=5$ having 2 irreducible components has type 3.
\end{prop}
\proof
The condition $t(C)=3$ is equivalent to $d_1+d_2=7$, which combined with \eqref{eqd-1} yields $d_1=3$, $d_2=d_3=4$.
Note that 
$$\tau(5,3)_{min}=4,$$
which implies $\tau(C) \geq 4$ with equality if and only if $C$ has type $3A$ with exponents $(3,4,4)$ by Theorem \ref{thmMT}. Consider now the case when $m=4$. Then
$$\tau(C)=18-2d_3-d_4 \text{ or } \tau(C)=18-d_3-2d_4$$
where $4=d_3 \leq d_4$. When $d_4=4$ we get the case $\tau(C)=6$ and $C$ of type $3B'$. When $d_4=5$, the first expression yields $\tau(C)=5$ and $C$ of type $3B$,
while the second expression gives $\tau(C)=4$. This last case is impossible by Theorem \ref{thmMT}.
Finally consider the case $m=5$. Then one has
$$\tau(C)=15-d_4-d_5,$$
with $4 \leq d_4 \leq d_5 \leq 6$. The last inequality here comes from \eqref{eqCM3}. Since $\tau(C) >4$ by Theorem \ref{thmMT}, we get the 4 possibilities for the exponents listed above when $C$ has type $3C$.

On the other hand, any nodal curve $C$ of degree $d=5$ having 2 irreducible components satisfies $d_1=3$ and $d_2=d_3=4$ by Theorem \ref{thmNodal} and \eqref{eqd-1}.
\endproof

\begin{cor}
\label{cord5}
A quintic curve $C$ with $\tau(C)=7$ is either a curve of type $2A$ with exponents $(3,3,4)$, or a curve of type $3C$ with exponents
$(3,4,4,4,4)$.
\end{cor}
\proof
Since 
$$\tau(5,4)_{max}'=4^2-\binom{5}{2}=6,$$
it follows that $d_1 \leq 3$. On the other hand
$$\tau(5,2)_{min}=8$$
and hence $d_1=3$. Recall that 
$$d_1d_2\geq (d-1)^2-\tau(C)=16-7=9,$$
see \cite[Theorem 2.2]{Euler}, with equality if and only if $m=3$ and $d_3=d-1$. Hence either the exponents are $(3,3,4)$ and then $C$ is a curve of type  $2A$, see \cite[Proposition 1.11]{ADP}, or $d_1=3$ and $d_2=4$. In this case the claim follows from
Proposition \ref{propd5}.
\endproof
\begin{ex}
\label{exd1+4}
(i) A nodal curve $C_{1,4}$ which is the intersection of a generic line with a smooth quartic has 4 nodes, and hence is a curve of type $3A$ by Proposition \ref{propd5}. 

(ii) A nodal curve $C'_{1,4}$ which is the intersection of a generic line with  quartic having one node, has 5 nodes in all, hence $\tau(C)=5$. In the special case when $C'_{1,4}$ is given by 
$$f=(x^4+y^4+xyz^2)z=0,$$
we get the exponents $(3,4,4,5)$, hence $C'_{1,4}$ has type $3B$.

(iii) A nodal curve $C''_{1,4}$ which is the intersection of a generic line with  quartic having 2 nodes, has 6 nodes in all, hence $\tau(C)=6$. In the special case when $C''_{1,4}$ is given by 
$$f=(xyz(x+y+z)+(x+y)^4)(x-y+z)=0,$$
we get the exponents $(3,4,4,4)$, hence $C''_{1,4}$ has type $3B'$.

(iv) A nodal curve $C'''_{1,4}$ which is the intersection of a generic line with  quartic having 3 nodes, has 7 nodes in all, hence $\tau(C)=7$
and hence $C'''_{1,4}$ has type $3C$ and the exponents are $(3,4,4,4,4)$.  As an example, one may take $C'''_{1,4}$  given by 
$$f=(x-y+z)(x^2y^2+y^2z^2+x^2z^2+xyz(x+y+z)) =0.$$
\end{ex}

\begin{ex}
\label{exd2+3}
(i) A nodal curve $C_{2,3}$ which is the intersection of a smooth conic and a smooth cubic has 6 nodes. The curve in Example \ref{exNodal2} (ii) is obtained in this way, has exponents $(3,4,4,4)$, and hence has type $3B'$.

(ii) A nodal curve $C'_{2,3}$ which is the intersection of a smooth conic with  a nodal cubic, has 7 nodes in all, hence $\tau(C)=7$. Therefore   $C''_{2,3}$ has type $3C$ and the exponents are $(3,4,4,4,4)$.

\end{ex}
The following example displays some non-nodal curves of type 3.
\begin{ex}
\label{exE7}

(i) By taking the quartic with one cusp $A_2$ from Example \ref{exd4} (i) and adding the line $x+y+z=0$, we get the quintic
$$C:f=(x+y+z)(x^4+x^2z^2+y^3z)=0.$$
This curve has 4 nodes, one cusp $A_2$ and hence $\tau(C)=6$.
A direct computation shows that the exponents are $(3,4,4,4)$, and hence $C$ has type $3B'$.

(ii) Here is an example of quintic curve of type $3C$ having an $E_7$ singularity.
Let
$$C:  f=x^5+y^5+yz(x^3+y^2z)=0.$$
Then a direct computation shows that $\tau(C)=7$, the exponents are $(3,4,4,4,4)$, and hence $C$ has type $3C$. The singularity at the point $(0:0:1)$ is clearly an $E_7$ singularity. Using Bézout Theorem, it is easy to see that this curve is irreducible.
\end{ex}

\begin{rk}
\label{rkd5}
(i) The quintics of type 3 constructed in Example \ref{exd1+4} (i), (ii), (iii) and (iv) and in Example \ref{exE7} (i) are obtained by adding a generic line to a quartic of type 3. It would be interesting to have a general result of this type.

(ii) We do not know if the cases $d=m=5$ and $5\leq \tau(C) \leq 6$ are possible.
\end{rk}

\section{On the line arrangements of type 3} 

For a line arrangement $\A:f=0$  with $d=\deg(f)$, it is known that
\begin{equation} 
\label{eqdm}
d_m \leq \reg D_0(f) \leq d-2,
\end{equation}
see \cite[Corollary 3.5]{Sch}. It follows that if $\A$ has type 3, then
$$d+2 =d_1+d_2 \leq 2d-4$$
and hence $d \geq 6$.  Moreover we know that
\begin{equation} 
\label{eqd1}
d_1 \leq d-m(\A),
\end{equation}
where $m(\A)$ denotes the maximal multiplicity of the singular points of $\A$, see \cite[Theorem 1.2]{Mich}.
Using these results we can prove the following.

\begin{prop}
\label{propA6}
An arrangement $\A$ of $d=6$ lines has type 3 if and only if it is a nodal arrangement. In this case $\tau(\A)=15$, $\A$ has type $3C$ and the exponents are
$(4,4,4,4,4)$.
\end{prop}
\proof
If $d=6$ it follows from \eqref{eqdm} that $d_1=d_2=4$. Then it follows from \eqref{eqd1} that $m(\A)=2$. We conclude using Example \ref{exC3} (ii).
\endproof
The case $d=7$ is also easy to treat.
\begin{prop}
\label{propA7}
An arrangement $\A$ of $d=7$ lines has type 3 if and only if it has only nodes except for exactly one triple point. In this case $\tau(\A)=22$, $\A$ has type $3C$ and the exponents are
$(4,5,5,5,5)$.
\end{prop}
\proof
If $d=7$ and $\A$ has type 3, then it follows from \eqref{eqdm} that $d_1=4$ and $d_2=5$.
Then \eqref{eqd1} implies that $\A$ has only double and triple points.
We use now \cite[Theorem 1.3]{Der} and conclude that there is exactly one triple point. Hence $\A$ has $n_2=18$ double points and $n_3=1$ triple points. It follows that
$$\tau(\A)=n_2+4n_3=22.$$
In view of \eqref{eqdm} the only possible exponents are $(4,5,5,5,5)$,
and hence $\A$ has type $3C$.
Conversely, assume that $n_2=18$ and $n_3=1$. Then \cite[Theorem 1.1]{Der} implies that
$$\dim D_0(f)_4=1$$
and this yields $d_1=4$ and $d_2=5$. Therefore $\A$ has type 3.
\endproof

\begin{ex}
\label{exA67}
(i) The line arrangement
$$\A:xyz(x+y+z)(x+2y+3z)(x+4y+9z)=0$$
is nodal and has $d=6$, hence by Proposition \ref{propA6} has type $3C$ and exponents $(4,4,4,4,4)$.

(ii) If we add the line $L:x-y=0$ to the arrangement $\A$ from the point (i) above, we get an arrangement $\A'$ having only double points except the triple point at $(0:0:1)$.  Hence by Proposition \ref{propA7} the arrangement $\A'$ has type $3C$, $\tau(\A')=22$ and exponents $(4,5,5,5,5)$.

\end{ex}

The case $d=8$ is definitely more subtle.

\begin{prop}
\label{propA8}
If an arrangement $\A$ of $d=8$ lines has type 3, then $\A$ satisfies one of the following three properties.

\begin{enumerate}
    \item $\A$ has only nodes except for exactly 4 triple points, $\tau(\A)=32$, $\A$ has type $3B'$ and the exponents are
$(5,5,5,5)$.

    \item $\A$ has only nodes except for exactly 3 triple points, $\tau(\A)=31$, $\A$ has type $3B$ and the exponents are
$(5,5,5,6)$.

    \item $\A$ has only nodes except for exactly 5  triple points, $\tau(\A)=33$, $\A$ has type $3C$ and the exponents are
$(5,5,5,5,5)$.

    \item $\A$ has only nodes except for exactly 4 triple points, $\tau(\A)=32$, $\A$ has type $3C$ and the exponents are
$(5,5,5,5,6)$.

    \item $\A$ has only nodes except for exactly 3  triple points, $\tau(\A)=31$, $\A$ has type $3C$ and the exponents are
$(5,5,5,6,6))$.

    \item $\A$ has only nodes except for exactly 2  triple points, $\tau(\A)=30$, $\A$ has type $3C$ and the exponents are
$(5,5,6,6,6))$.

 \item $\A$ has only nodes except for exactly one point of multiplicity 4,  $\tau(\A)=31$, $\A$ has type $3C$ and the exponents are
$(4,6,6,6,6)$.

\end{enumerate}

\end{prop}

\proof
If $d=8$, then it follows from \eqref{eqdm} that either $d_1=d_2=5$ or $d_1=4$ and $d_2=6$.
Consider first the case $d_1=d_2=5$. Then it follows from 
\eqref{eqdm} that $\A$ has as exponents
$(5,5,5,5)$, $(5,5,5,6)$,  $(5,5,5,5,5)$,  $(5,5,5,5,6)$, $(5,5,5,6,6)$ or $(5,5,6,6,6)$.

In each case $\A$ has $d_1=5=8-3$, and hence there are no points of multiplicity $>3$ by \eqref{eqd1}. Moreover, the number of triple points
$n_3$ is given by the number of 5's in the corresponding exponent  by  \cite[Theorem 1.3]{Der}. Then formula number of nodes is given by the well know
$$n_2=\binom{8}{2}-n_3 \binom{3}{2}.$$
It follows that 
$$\tau(\A)=n_2+4n_3$$
which coincides with the values given by Theorem \ref{thm1}.

Consider now the case $d_1=4$ and $d_2=6$. Then by \eqref{eqdm} the only possible exponents are $(4,6,6,6,6)$ and $\A$ has only double points except for one point of multiplicity 4 by  \cite[Theorem 1.3]{Der}.
Then $\A$ has 
$$n_2=28-\binom{4}{2}=22$$
double points and
$$\tau(\A)=n_2+9n_4=31,$$
which coincides with the value given by Theorem \ref{thm1}, (iv).
\endproof

\begin{ex}
\label{exA8}

(i) The line arrangement
$$\A:f=xy(x-y-z)(x-y+z)(2x+y-2z)(x+3y-3z)(3x+2y+3z)(x+5y+5z)=0$$
obtained from the second line arrangement in Example \ref{exZ} by deleting the last factor has $d=8$, $n_2=16$, $n_3=4$ and exponents
$(5,5,5,5)$ as in Proposition \ref{propA8} (i).

(ii) The line arrangement
$$\A:f=(x+y)(x+2y)(x+3y)(x-z)(x-2z)(x-3z)(9x-2y-z)(10x-y-z)=0$$
has 
$d=8$, $n_2=19$, $n_3=3$, $\tau(\A)=31$ and exponents
$(5,5,5,6)$ as in Proposition \ref{propA8} (2).

(iii) The line arrangement
$$\A: f=xyz(x-z)(y-z)(x-2y)(y-z+3x)(x-z+7y)=0$$
has 
$d=8$, $n_2=13$, $n_3=5$, $\tau(\A)=33$ and exponents
$(5,5,5,5,5)$ as in Proposition \ref{propA8} (3).

(iv)The line arrangement
$$\A: f=xyz(x-z)(y-z)(x-2y)(2x-y-z)(x+7y+11z)=0$$
has 
$d=8$, $n_2=16$, $n_3=4$, $\tau(\A)=32$ and exponents
$(5,5,5,5,6)$ as in Proposition \ref{propA8} (4).

(v) The line arrangement
$$\A: f=x(x-z)y(y-z)(x+3y)z(x+2y+3z)(x+4y+9z)=0$$
 has $d=8$, $n_2=19$, $n_3=3$, $\tau(\A)=31$ and exponents
$(5,5,5,6,6)$ as in Proposition \ref{propA8} (5).

(vi) The line arrangement
$$\A: f=yz(x^2-y^2)(x^2-z^2)(x-2y+13z)(x+7y+11z)=0$$
 has $d=8$, $n_2=22$, $n_3=2$, $\tau(\A)=30$ and exponents
$(5,5,6,6,6)$ as in Proposition \ref{propA8} (6).
Here the triple points are $(0:0:1)$ and $(0:1:0)$ and they are not connected by a line in $\A$. On the other hand, the arrangement
$$\A':f=x(x+y)(x-y)(x^2-z^2)(x-2y+13z)(x+7y+11z)(y+13z)=0$$
has the same numerical invariants as $\A$, the same triple points, but this time they are connected by the line $x=0$ which is in $\A'$.

(vii) The line arrangement
$$\A: f=  xyz(x+y+z)(x+2y+3z)(x+4y+9z)(x-2y)(x-3y)=0$$
has
$d=8$, $n_2=22$, $n_4=1$, $\tau(\A)=31$ and exponents
$(4,6,6,6,6)$ as in Proposition \ref{propA8} (7).

\end{ex}

\begin{rk}
\label{rkA8}
(i)  The line arrangement
$$\B:xyz(x-y)(y-z)(x-z)=0$$
has $d=6$, $n_2=3$ double points and $n_3=4$ triple points.
To get an arrangement $\B'$ having $d=8$, $n_2=16$ and $n_3=4$ we may add to $\B$ two generic lines, for instance  $x+y+z=0$ and $x+3y+7z=0$. This arrangement $\B'$ has $\tau(\B')=32$ but the exponents are $(4,5,6,6)$ as a direct computation with SINGULAR shows. Hence the arrangement $\B'$ has type $2B$ and not type 3, even if $\B'$ has the same pair of numerical invariants $(n_2,n_3)$,
sometimes called the {\it weak combinatorial type}, as the arrangement $\A$ from Example \ref{exA8}
(i) above.

If we add to $\B$ the lines $x+y+z=0$ and $x+2y+3z=0$ we get a new line arrangement $\B''$ having $n_2=13 $, $n_3=5$, $\tau(\B'')=33$ and exponents $(4,5,5)$. Therefore $\B''$ has type $2A$ and there is no converse claim in
Proposition \ref{propA8} (3), even if the Tjurina number 33 occurs only once in the cases listed in this Proposition.

(ii) We see in Proposition \ref{propA8} that the weak combinatorial type does not determine completely the exponents even inside type 3 line arrangements of degree 8. The arrangements in (1) and (4) (resp. (2) and (5))  have the same weak combinatorial type, but distinct intersection lattices. Unfortunately, it seems that  there is no general approach allowing one to get the exponents from the intersection lattice, even in the cases when this lattice determines the exponents.

\end{rk}
The only converse result related to Proposition \ref{propA8} which we can prove is the following. It shows that the weak combinatorial type and the exponents determine each other in the situation of Proposition \ref{propA8} (6),
even if we have seen in Example \ref{exA8} (vi) that the intersection lattice is not determined by the exponents.

\begin{cor}
\label{corA8}
If an arrangement $\A$ of $8$ lines has only nodes except exactly 2 triple points, then $\A$ has type $3C$ and exponents $(5,5,6,6,6)$.
\end{cor}
\proof
In view of Proposition \ref{propA8} it is enough to show that $d_1=d_2=5$.
The hypothesis on $\A$ and  \cite[Theorem 1.3]{Der} imply that
$$\dim D_0(f)_5=2.$$
This implies that $d_1=d_2=5$. Indeed, if $d_1<5$, then
$$\dim D_0(f)_5 \geq \dim S_1=3,$$
again by \cite[Theorem 1.3]{Der}. Hence $d_1=5$, and $d_2=5$ since we need a second generators in $D_0(f)_5$.
\endproof

Note that \eqref{eqdm} implies that
$d_2 \leq d-2$, and hence for a type 3 line arrangement $\A$ one has
$$d_1=d+2-d_2 \geq 4.$$
The following result shows that this is the only condition on the first two exponents of a type 3 line arrangement.
\begin{thm}
\label{thmAd}
Let $\B_1$ (resp. $\B_2$) be the arrangement consisting  of $n_1$ lines passing through a point $a_1 \in \PP^2$ (resp. $n_2$ lines  passing through a point $a_2 \in \PP^2 \setminus \B_1$) and such that $2 \leq n_1 \leq n_2$. Then the arrangement $\A:f=0$ obtained from the arrangement $\B=\B_1 \cup \B_2$  by adding two generic lines has type $3C$ and  exponents 
$$(n_1+2,n_2+2,n_1+n_2, n_1+n_2,n_1+n_2).$$ 
\end{thm}

\proof
Let $\B':f'=0$ be the arrangement  obtained from the arrangement $\B$  by adding one generic line. Then \cite[Theorem 1.6]{ADP} shows that 
$\B'$ has exponents as first exponents $d_1'=n_1+1$ and $d_2'=n_2+1$. Indeed, the exponents of $\B$ are $(n_1,n_2,n_1+n_2-2)$ by \cite[Proposition 1.17]{ADP} and hence the condition
$$n_2 \leq \deg \B -2=n_1+n_2-2$$
in \cite[Theorem 1.6]{ADP} is fulfilled. It follows that $\B'$ has type 2,
and in fact for $n_1 \geq 3$, $\B'$ has type $2B$ and exponents
$$(n_1+1,n_2+1,n_1+n_2-1, n_1+n_2-1),$$
see \cite[Corollary 1.18]{ADP}. For $n_1=2$, using the formulas for $\tau(\B')$ given in \cite[Proposition 1.11]{ADP} and the obvious fact that
$$\tau(\B')=n_2^2+n_2+4$$
in this case, we get that $\B'$ has type $2B$ and the exponents are given by the above formulas.

We apply now the exact sequence
\cite[(2.6)]{ADP} to the arrangement $\B'$ and a new generic line $L$. Then the intersection $\B' \cap L$ consists of $r=\deg \B'=n_1+n_2+1$ points and we get in this way an exact sequence
$$0 \to D_0(f')_{k-1} \to D_0(f)_k \to H^0(L,\OO_L(k+1-r))$$
for any integer $k$.
Since 
$$H^0(L,\OO_L(k+1-r)=H^0(L,\OO_L(k-n_1-n_2))=0$$
for $k<n_1+n_2$ it follows that the first 2 exponents of $\A$ are
$d_1=d_1'+1=n_1+2$ and $d_2=d_2'+1=n_2+2$. The other two generators of $D_0(f')$ gives rise to  elements of degree $n_1+n_2$, and similar elements can be added by $H^0(L,\OO_L(k-n_1-n_2))$ when $k \geq n_1+n_2$. Since $\A$ contains $d=n_1+n_2+2$ lines, it follows that all generators have degree $\leq n_1+n_2$.
Therefore $\A$ has type $3$ and the exponents are either
$$(n_1+2,n_2+2,n_1+n_2, n_1+n_2) \text{ or } (n_1+2,n_2+2,n_1+n_2, n_1+n_2, n_1+n_2).$$
A direct computation, using the obvious facts that 
$$\tau(\B)=(n_1-1)^2+(n_2-1)^2+n_1n_2, \ \tau(\B')=\tau(\B)+n_1+n_2 \text{ and } \tau(\A)=\tau(\B')+n_1+n_2+1$$
plus the formulas in Theorem \ref{thm1} shows that only the second case is possible. This completes the proof of our claim.
\endproof

\end{document}